\documentclass[11pt,a4paper]{article}

\usepackage{amsmath,amsfonts,amssymb,latexsym,graphics,epsfig,url}
\usepackage{color}
\usepackage{amsthm}
\usepackage[english]{babel}
\usepackage{graphicx}
\usepackage{soul}
\usepackage[utf8]{inputenc}

\newtheorem{theorem}{Theorem}[section]

\newtheorem{lemma}[theorem]{Lemma}
\newtheorem{corollary}[theorem]{Corollary}

\newtheorem{conjecture}[theorem]{Conjecture}
\newtheorem{example}[theorem]{Example}

\DeclareMathOperator{\dist}{dist}
\DeclareMathOperator{\spec}{sp}

\def\u{\mbox{\boldmath $u$}}
\def\v{\mbox{\boldmath $v$}}

\def\w{\mbox{\boldmath $w$}}

\def\vec0{\mbox{\boldmath $0$}}

\def\B{\mbox{\boldmath $B$}}

\def\I{\mbox{\boldmath $I$}}
\def\J{\mbox{\boldmath $J$}}

\def\L{\mbox{\boldmath $L$}}
\def\M{\mbox{\boldmath $M$}}
\def\O{\mbox{\boldmath $O$}}

\def\I{\mbox{\boldmath $I$}}
\def\J{\mbox{\boldmath $J$}}

\def\1{\mbox{\boldmath $1$}}

\def\Re{\mathbb R}

\newcommand\restr[2]{\ensuremath{\left.#1\right|_{#2}}}

\begin{document}
	
	\title{On the algebraic connectivity of token graphs
		\thanks{This research of C. Dalf\'o and M. A. Fiol has been partially supported by
			AGAUR from the Catalan Government under project 2017SGR1087 and by MICINN from the Spanish Government under project PGC2018-095471-B-I00. The research of C. Dalf\'o has also been supported by MICINN from the Spanish Government under project MTM2017-83271-R.}
	}
	\author{C. Dalf\'o$^a$,  M. A. Fiol$^b$,\\
		\\
		{\small $^a$Dept. de Matem\`atica, Universitat de Lleida, Igualada (Barcelona), Catalonia}\\
		{\small {\tt cristina.dalfo@udl.cat}}\\
		{\small $^{b}$Dept. de Matem\`atiques, Universitat Polit\`ecnica de Catalunya, Barcelona, Catalonia} \\
		{\small Barcelona Graduate School of Mathematics} \\
		{\small  Institut de Matem\`atiques de la UPC-BarcelonaTech (IMTech)}\\
		{\small {\tt miguel.angel.fiol@upc.edu} }\\
	}

	\date{}
	\maketitle
	
	\begin{abstract}
		We study the algebraic connectivity (or second Laplacian eigenvalue) of token graphs, also called symmetric powers of graphs. The $k$-token graph $F_k(G)$ of a graph $G$ is the
		graph whose vertices are the $k$-subsets of vertices from $G$, two of which being adjacent whenever their symmetric difference is a pair of adjacent vertices in $G$.
		Recently, it was conjectured that the algebraic connectivity of $F_k(G)$ equals the algebraic connectivity of $G$. In this paper, we prove the conjecture for new infinite families of graphs, such as trees and graphs with  maximum degree large enough.
	\end{abstract}
	
	
	
	\section{Introduction}
	\label{sec:-1}
\begin{figure}[t]
		\begin{center}
			\includegraphics[width=5cm]{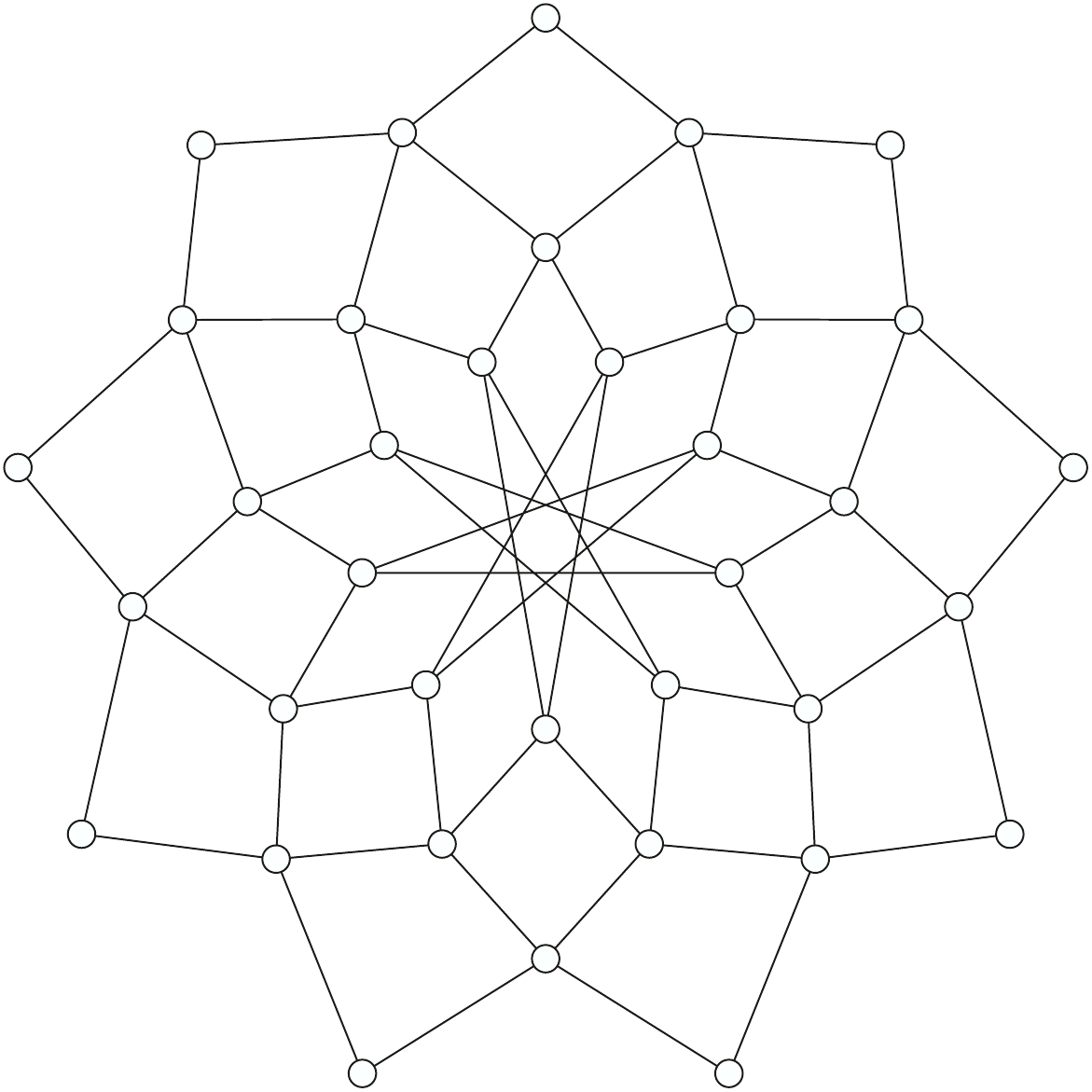}
			\caption{The $2$-token graph $F_2(C_9)$ of the cycle graph, with vertex set $V(C_9)=\{0,1,\ldots,8\}$.
				The vertices on the circumference of radius  $r_{\ell}$, with $\ell=1,2,3,4$ and $r_1>r_2>r_3>r_4$  are $\{i,j\}$ with $\dist(i,j)=\ell$ in $C_9$.}
			\label{fig1}
		\end{center}
	\end{figure}

	Let $G$ be a simple graph with vertex set $V(G)=\{1,2,\ldots,n\}$ and edge set $E(G)$. Let $\Delta(G)$ denote the maximum degree of $G$. For a given integer $k$ such that
	$1\le k \le n$, the {\em $k$-token graph} $F_k(G)$ of $G$ is the graph whose vertex set $V (F_k(G))$ consists of the ${n \choose k}$
	$k$-subsets of vertices of $G$, and two vertices $A$ and $B$
	of $F_k(G)$ are adjacent whenever their symmetric difference $A \bigtriangleup B$ is a pair $\{a,b\}$ such that $a\in A$, $b\in B$, and $(a,b)\in E(G)$.
	The naming `token graph'
	comes from an observation in
	Fabila-Monroy,  Flores-Pe\~{n}aloza,  Huemer,  Hurtado,  Urrutia, and  Wood \cite{ffhhuw12}, that vertices of $F_k(G)$ correspond to configurations
	of $k$ indistinguishable tokens placed at distinct vertices of $G$, where
	two configurations are adjacent whenever one configuration can be reached
	from the other by moving one token along an edge from its current position
	to an unoccupied vertex. Thus,
	the maximum degree of $F_k(G)$ satisfies
	\begin{equation}
		\label{DeltaFk}
		\Delta(F_k(G))\le k\Delta(G).
	\end{equation}
	In Figure \ref{fig1}, we show the 2-token graph of the cycle $C_9$ on 9 vertices.
	
	Note that if $k=1$, then $F_1(G)\cong G$; and if $G$ is the complete graph $K_n$, then $F_k(K_n)\cong J(n,k)$, where $J(n,k)$ denotes the Johnson graph~\cite{ffhhuw12}.
	
	Token graphs have some applications in physics. For instance, a
	relationship between token graphs and the exchange of Hamiltonian operators in
	quantum mechanics is given in Audenaert, Godsil, Royle, and Rudolph \cite{agrr07}.
	
	Recently, it was conjectured by Dalf\'o, Duque, Fabila-Monroy, Fiol, Huemer, Trujillo-Negrete, and Zaragoza Mart\'{\i}nez \cite{ddffhtz21} that the algebraic connectivity of $F_k(G)$ equals the algebraic connectivity of $G$. In this paper, we prove the conjecture for new infinite families of graphs, such as trees, multipartite complete graphs, and graphs with large enough maximum degree.
	
	\section{Known results}
	Let us first introduce some notation and known results that are used throughout the paper.
	The transpose of a matrix $\M$ is denoted by $\M^\top$, the
	identity matrix by $\I$, the all-$1$ vector $(1,..., 1)^{\top}$ by $\1$, the all-$1$ (universal) matrix  by $\J$, and the all-$0$ vector and all-$0$ matrix by $\vec0$
	and $\O$, respectively.
	Let $[n]:=\{1,\ldots,n\}$ and ${[n]\choose k}$ denote the set of $k$-subsets of $[n]$, which is the set of vertices of the $k$-token graph.
	
	For our purpose, it is convenient to denote by $W_n$ the set of all column vectors $\v$ such that  $\v^{\top }\1 = 0$.
	Recall that any square matrix $\M$ with all zero row sums has an eigenvalue $0$ with corresponding eigenvector $\1$.
	
	\subsection{The algebraic connectivity of token graphs}
	
	When $\M=\L(G)$, the Laplacian matrix of a graph $G$, the matrix is positive semidefinite, with eigenvalues $(0=)\lambda_1\le \lambda_2\le \cdots \le \lambda_n$. Its second smallest eigenvalue $\lambda_2$ is known as the {\em algebraic connectivity} of  $G$ (see Fiedler \cite{fi73}), and we denote it by $\alpha(G)$.
	The spectral radius $\lambda_{\max}(G)=\lambda_n$ satisfies several lower and upper bounds (see
	Patra and  Sahoo \cite{ps17} for a survey). Here, we will use the following ones in terms of the maximum degree of $G$:
	\begin{equation}
		\label{bound-sp-rad}
		1+\Delta(G)\le \lambda_{\max}(G)\le 2\Delta(G).
	\end{equation}
	The upper bound is due to Fiedler \cite{fi73}, whereas the lower bound was proved by Grone and  Merris in \cite{gm94}.
	
	In this paper, we want to study the algebraic connectivity of token graphs. As far as we know, this study was initiated by Dalf\'o,  Duque, Fabila-Monroy,  Fiol, Huemer,  Trujillo-Negrete, and  Zaragoza Mart\'{\i}nez in \cite{ddffhtz21}, where they proved the following result.
	
	\begin{lemma}[\cite{ddffhtz21}]
		\label{coro:LkL1}
		Let $G$ be a graph with Laplacian matrix $L_1$. Let $F_k=F_k(G)$ be its token graph with Laplacian $L_k$.
		Let $\B$ be the so-called $(n;k)$-\emph{binomial matrix}, which is an ${n \choose k}\times n$ matrix whose rows are the characteristic vectors of the $k$-subsets of $[n]$ in a given order.
		Then, the following holds:
		\begin{itemize}
			\item[$(i)$]
			If $\v$ is a $\lambda$-eigenvector of $\L_1$, then $\B\v$ is a $\lambda$-eigenvector of $\L_k$.
			Thus, the Laplacian spectrum (eigenvalues and their multiplicities) of $\L_1$ is contained in the Laplacian spectrum of $\L_k$.
			\item[$(ii)$]
			If $\u$ is a $\lambda$-eigenvector of $\L_k$ such that $\B^{\top}\u\neq \vec0$, then $\B^{\top}\u$
			is a $\lambda$-eigenvector of $\L_1$.
		\end{itemize}
	\end{lemma}
	
	Given two integers $n,k$ such that $k\in [n]$, the {\em Johnson graph} $J(n,k)$ can be  defined as the $k$-token graph of the complete graph $K_n$ , $F_k(K_n)\cong J(n,k)$. It is known that these graphs are antipodal (but not bipartite) distance-regular graphs, with degree $d=k(n-k)$, diameter $D=\min\{k,n-k\}$, and Laplacian spectrum
	(eigenvalues and multiplicities)
	\begin{equation}
		\label{spJ(n,k)}
		\lambda_j=d-\mu_j=j(n+1-j)\qquad \mbox{and} \qquad m_j={n\choose j}-{n\choose j-1},\qquad j=0,1,\ldots,D.
	\end{equation}
	(See again \cite{ddffhtz21}).
	For example, $F_2(K_{4})\cong J(4,2)$ is a $2$-regular graph with $n=6$ vertices, diameter $D=2$, and Laplacian spectrum
	$
	S( F_2(K_{4}))=\{0^1, 4^{3}, 6^{2}\}$.
	
	Let us consider a graph $G$ and its complement $\overline{G}$, with respective Laplacian matrices $L_G$ and $L_{\overline{G}}$.
	Since $L_G+L_{\overline{G}}=n\I-\J$, the Laplacian spectrum of $\overline{G}$ is the complement of the Laplacian spectrum of $G$ with respect to the Laplacian spectrum of the complete graph $K_n$. We represent this as
	$$
	\spec G \oplus \spec \overline{G} = \spec K_n,
	$$
	where each eigenvalue of $G$ and each eigenvalue of $\overline{G}$ are used once.
	In \cite{ddffhtz21}, is was shown that a similar relationship holds between the Laplacian spectra of the $k$-token of $G$ and the $k$-token of $\overline{G}$, but now with respect to the Laplacian spectrum of the Johnson graph.
	\begin{theorem}[\cite{ddffhtz21}]
		\label{theo:pairing}
		Let $G=(V,E)$ be a graph on $n=|V|$ vertices, and let $\overline{G}$ be its complement. For a given $k$, with $1\leq k\le n-1$, let us consider the token graphs $F_k(G)$ and  $F_k(\overline{G})$. Then, the Laplacian spectrum of $F_k(\overline{G})$ is the complement of the Laplacian spectrum of $F_k(G)$ with respect to the Laplacian spectrum of the Johnson graph $J(n,k)=F_k(K_n)$.
		That is, every eigenvalue  $\lambda_J$ of $J(n,k)$ is the sum of one eigenvalue $\lambda_{F_k(G)}$ of $F_k(G)$ and one eigenvalue $\lambda_{F_k(\overline{G})}$ of $F_k(\overline{G})$, where each $\lambda_{F_k({G})}$ and each $\lambda_{F_k(\overline{G})}$ is used once:
		\begin{equation}
			\label{spFk(G)-sp(Fk(noG))}
			\lambda_{F_k({G})}+\lambda_{F_k(\overline{G})}=\lambda_J.
		\end{equation}
	\end{theorem}
	
	As a consequence of Lemma \ref{coro:LkL1}, the spectrum of $J(n,k)$ in \eqref{spJ(n,k)}, and Theorem \ref{theo:pairing}, we can state the following lemma where, for simplicity, we assume that both $G$ and $\overline{G}$ are connected.
	\begin{lemma}
		\label{partition}
		Let $\Lambda$ be the set of pairs $(\lambda,\overline{\lambda})$ of eigenvalues of $F_k(G)$ and $F_k(\overline{G})$, with $k\le n/2$, sharing both the same eigenvector $\v$ with $J(n,k)$. Then $\Lambda$ can be partitioned into the sets $\Lambda_0,\Lambda_1,\ldots,\Lambda_k$ such that
		$\Lambda_0=\{(0,0)\}$, and
		$\Lambda_j=\{(\lambda,\overline{\lambda}):\lambda+\overline{\lambda}=j(n+1-j)\}$ for $j=1,\ldots,k$.
		Moreover, the eigenvectors and eigenvalues of each set satisfy: $\v=\1$ in $\Lambda_0$; the eigenvalues in $\Lambda_1$ correspond to the eigenvalues of $G=F_1(G)$ and $\overline{G}=F_1(\overline{G})$; the eigenvalues in $\Lambda_j$ come from $\spec F_j(G)\setminus \spec F_{j-1}(G)$ and $\spec F_j(\overline{G})\setminus \spec F_{j-1}(\overline{G})$ with eigenvectors $\v$  such that $\B^{\top}\v=\vec0$, for $j=2,\ldots, k$.
	\end{lemma}
	
	Let us show an example of the results in Theorem \ref{theo:pairing} and Lemma \ref{partition}.
	
	\begin{example}
		\label{exemple}
		Consider the graph $G$ and its complement graph $\overline{G}$ of Figure \ref{fig2}. 
		\begin{figure}[h]
			\begin{center}
				\includegraphics[width=7cm]{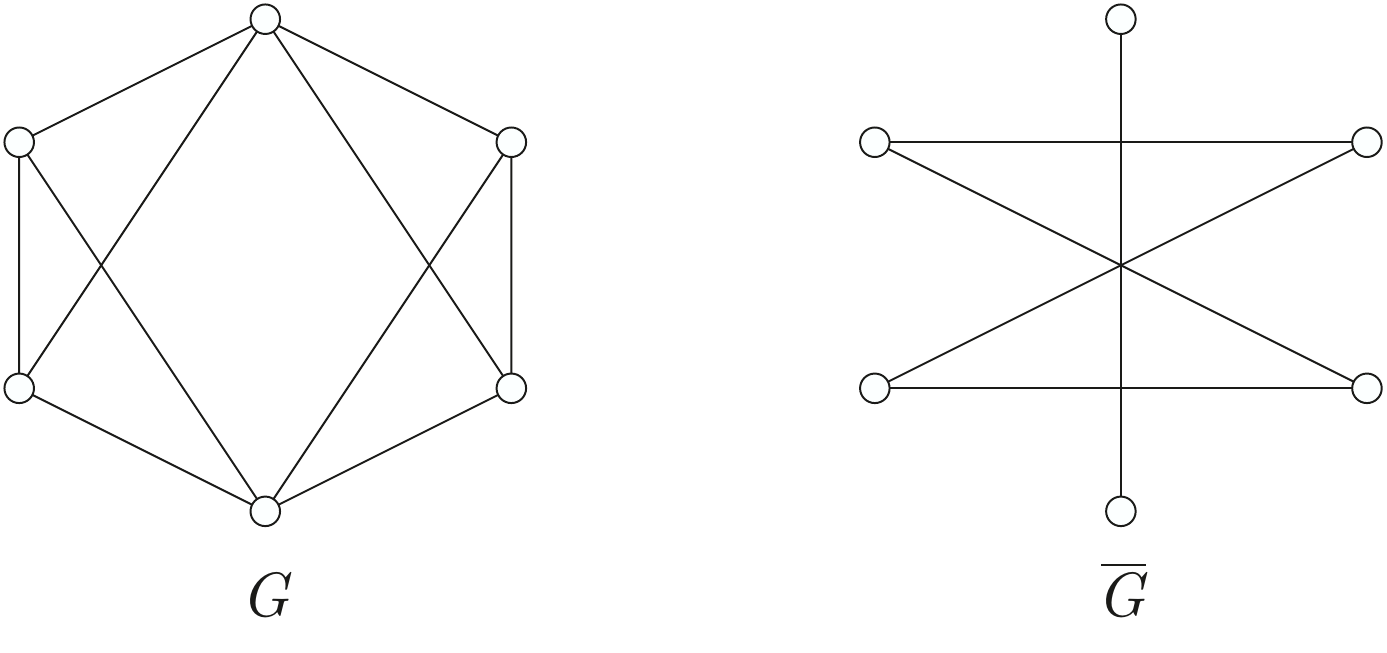}
				\caption{The graph $G$ and its complement graph $\overline{G}$ of Example \ref{exemple}.}
				\label{fig2}
			\end{center}
		\end{figure}
		The spectra of $G$, $\overline{G}$, and their $k$-tokens for $k=2,3$ are the following:
		\begin{align*}
			\spec G &=\{0,2,4^{[3]},6\}\subset \spec F_2(G) =\{0,2,4^{[5]},6^{[4]},8^{[3]},10\}\\
			&\subset \spec F_3(G) =\{0,2,4^{[6]},6^{[4]},8^{[5]},10^{[3]}\}.\\
			\spec \overline{G} &=\{0^{[2]},2^{[3]},4\}\subset \spec F_2(\overline{G})=\{0^{[3]},2^{[6]},4^{[4]},6^{[2]}\} \subset
			\spec F_3(\overline{G}) =\{0^{[3]},2^{[8]},4^{[6]},6^{[2]}\}.
		\end{align*}
		Then, as shown in Table \ref{table1}, there is a pairing between the eigenvalues of $F_3(G)$ and the eigenvalues of $F_3(\overline{G})$ satisfying Theorem \ref{theo:pairing}.
		Namely, 
		$$
		\spec F_3(G)\oplus \spec F_3(\overline{G}) = \spec J(6,3)=\{0, 6^{[5]}, 10^{[9]}, 12^{[5]}\}.
		$$
		Thus, the pairs of $\Lambda_0$, $\Lambda_1$, $\Lambda_2$, and $\Lambda_3$ add up to 0, 6, 10, and 12, respectively.
		\begin{table}
			\begin{center}
				\scriptsize
				\setlength\tabcolsep{3pt}
				\begin{tabular}{|c||c|c|c|c|c|c|c|c|c|c|c|c|c|c|c|c|c|c|c|c|}
					\hline
					$F_3(\overline{G})\setminus F_3(G)$
					& 0 & 2 & 4 & 4 & 4 & $\mathit{4}$ & $\mathit{4}$ & $\mathbf{4}$ & $\mathit{6}$ & $\mathit{6}$ & $\mathit{6}$ & $\mathit{6}$ & $\mathit{8}$ & $\mathit{8}$ & $\mathit{8}$ & $\mathbf{8}$ & $\mathbf{8}$ & $\mathit{10}$ & $\mathbf{10}$ & $\mathbf{10}$\\
					\hline\hline
					0
					& 0 &   &   &   &   &   &   &   &   &   &   &   &   &   &   &   &   &    &    &   \\
					\hline
					0
					&   &   &   &   &   &   &   &   & 6 &   &   &   &   &   &   &   &   &    &    &   \\
					\hline
					$\mathit{0}$
					&   &   &   &   &   &   &   &   &   &   &   &   &   &   &   &   &   & $\mathit{10}$   &    &   \\
					\hline
					2
					&   &   & 6  &  &   &   &   &   &   &   &   &   &    &    &   &   &   &    &    &   \\
					\hline
					2
					&   &   &   & 6 &  &   &   &   &   &   &   &   &    &    &   &   &   &    &    &   \\
					\hline
					2
					&   &   &   &   & 6 &   &   &   &   &   &   &   &   &    &    &   &   &    &    &   \\
					\hline
					$\mathit{2}$
					&   &   &   &   &   &   &   &   &   &   &   &   & $\mathit{10}$   &  &    &   &   &    &    &   \\
					\hline
					$\mathit{2}$
					&   &   &   &   &   &   &   &   &   &   &   &   &    &  $\mathit{10}$  &   &   &   &    &    &   \\
					\hline
					$\mathit{2}$
					&   &   &   &   &   &   &   &   &   &   &   &   &   &   & $\mathit{10}$  &   &   &    &    &   \\
					\hline
					$\mathbf{2}$
					&   &   &   &   &   &   &   &   &   &   &   &   &   &   &   &   &   &    &  $\mathbf{12}$  &  \\
					\hline
					$\mathbf{2}$
					&   &   &   &   &   &   &   &   &   &   &   &   &   &   &   &   &   &    &    &  $\mathbf{12}$ \\
					\hline
					4
					&   & 6  &   &   &   &   &   &   &   &   &   &   &   &   &   &   &   &    &    &   \\
					\hline
					$\mathit{4}$
					&   &   &   &   &   &   &   &   &   & $\mathit{10}$ &   &   &   &   &   &   &   &    &    &   \\
					\hline
					$\mathit{4}$
					&   &   &   &   &   &   &   &   &   &   & $\mathit{10}$ &   &   &   &   &   &   &    &    &   \\
					\hline
					$\mathit{4}$
					&   &   &   &   &   &   &   &   &   &   &   &  $\mathit{10}$ &   &   &   &   &   &    &    &   \\
					\hline
					$\mathbf{4}$
					&   &   &   &   &   &   &   &   &   &   &   &   &   &   &   & $\mathbf{12}$  &   &    &    &   \\
					\hline
					$\mathbf{4}$
					&   &   &   &   &   &   &   &   &   &   &   &   &   &   &   &   & $\mathbf{12}$  &    &    &   \\
					\hline
					$\mathit{6}$
					&   &   &   &   &   &  $\mathit{10}$ &   &   &   &   &   &   &   &   &   &   &   &    &    &   \\
					\hline
					$\mathit{6}$
					&   &   &   &   &   &   & $\mathit{10}$  &   &   &   &   &   &   &   &   &   &   &    &    &   \\
					\hline
					$\mathbf{8}$
					&   &   &   &   &   &   &   &  $\mathbf{12}$ &   &   &   &   &   &   &   &   &   &    &    &   \\
					\hline
				\end{tabular}
			\end{center}
			\caption{The spectra of $F_3(G)$ and $F_3(\overline{G})$ giving the spectrum of $J(6,3)$ by addition. The eigenvalues of $\Lambda_2$ are in italics. The eigenvalues of $\Lambda_3$ are in boldface. The corresponding additions, giving the eigenvalues of $J(6,3)$, are written accordingly.}
			\label{table1}
		\end{table}
		
	\end{example}
	
	Concerning the algebraic connectivity of token graphs,
	the authors in \cite{ddffhtz21} proposed the following conjecture.
	
	\begin{conjecture}[\cite{ddffhtz21}]
		\label{conjecture}
		Let $G$ be a graph on $n$ vertices. Then, for every $k=1,\ldots,n-1$, the algebraic connectivity of its token graph $F_k(G)$ equals the algebraic connectivity of $G$.
	\end{conjecture}
	
	As a consequence of Lemma \ref{coro:LkL1}$(i)$, and that $F_k(G)=F_{n-k}(G)$,  the conjecture only needs to be proved for the case $k=\lfloor n/2 \rfloor$.
	Moreover,  it was noted that the conjecture
	also holds when the graph $G$ is disconnected and for those graphs whose token graphs are regular, which are $K_n$, $S_n$ (with even $n$ and $k=n/2$), and their complements.
	Also, computer exploration showed that $\alpha(F_2(G))=\alpha(G)$ for all graphs with at most $8$ vertices.
	Moreover, it was shown that the conjecture also holds for the following infinite families of graphs.
	\begin{theorem}[\cite{ddffhtz21}]
		\label{theo:alg-connec-antic}
		For each  of the following classes of graphs, the algebraic connectivity of a token graph $F_k(G)$ equals the algebraic connectivity of $G$.
		\begin{itemize}
			\item[$(i)$]
			Let $G=K_n$ be the complete graph on $n$ vertices. Then,
			$\alpha(F_k(G))=\alpha(G)=n$ for every $n$ and $k=1,\ldots,n-1$.
			\item[$(ii)$]
			Let $G=S_n$ be the star graph on $n$ vertices. Then,
			$\alpha(F_k(G))=\alpha(G)=1$ for every $n$ and $k=1,\ldots,n-1$.
			\item[$(iii)$]
			Let $G=P_n$ be the path graph on $n$ vertices. Then, $\alpha(F_k(G))=\alpha(G)\linebreak =2(1-\cos(\pi/n))$ for every $n$ and  $k=1,\ldots, n-1$.
			\item[$(iv)$]
			Let $G= K_{n_1,n_2}$ be the complete bipartite graph on $n=n_1+n_2$ vertices, with $n_1\le n_2$. Then, $\alpha(F_k(G))=\alpha(G)=n_1$ for every $n_1,n_2$ and $k=1,\ldots,n-1$.
		\end{itemize}
	\end{theorem}

	\section{New results}
	\label{sec:-7}
	
	In this section, we prove more results supporting Conjecture \ref{conjecture}. In our proofs, we use the following concepts and lemma.
	
	Given a graph $G=(V,E)$ of order $n$,
	we say that a vector $\v\in \mathbb{R}^n$ is an \textit{embedding} of $G$ if $\v\in W_n$.
	Note that if $\v$ is a $\lambda$-eigenvector of $G$, with $\lambda>0$, then it is an embedding of $G$.
	
	For a graph $G$ with Laplacian matrix $\L(G)$, and an embedding $\v$ of $G$, let
	\begin{equation}
		\label{rayleigh-quotient}
		\lambda_G(\v):=\frac{\v^{\top}\L(G)\v}{\v^{\top}\v}=\frac{\sum\limits_{(i,j)\in E}[\v(i)-\v(j)]^2}{\sum\limits_{i\in V}\v^2(i)}
	\end{equation}
	where
	$\v(i)$ denotes the entry of $\v$ corresponding to the vertex $i\in V(G)$.
	The value of $\lambda_G(\v)$ is known as
	the {\em Rayleigh quotient}.
	If $\v$ is an eigenvector of $G$, then its corresponding eigenvalue is $\lambda(\v)$.
	Moreover, for an embedding $\v$ of $G$, we have
	\begin{equation}
		\label{bound-lambda(v)}
		\alpha(G)\le \lambda_G(\v),
	\end{equation}
	and we have equality  when $\v$ is an $\alpha(G)$-eigenvector of $G$.
	\\
	\begin{lemma}
		\label{lem:-vertex}
		Let $G^+=(V^+,E^+)$ be a graph on the vertex set $V=\{1,2,\ldots,n+1\}$, having a vertex of degree $1$, say the vertex $n+1$ that is adjacent to $n$.
		Let $G=(V,E)$ be the graph obtained from $G^+$ by deleting the vertex $n+1$.
		Then,
		$$
		\alpha(G)\ge \alpha(G^+),
		$$
		with equality if and only if the $\alpha(G)$-eigenvector $\v$ of $G$ has entry $\v(n)=0$.
	\end{lemma}
	\begin{proof}
		Let $\v\in W_n$ be an eigenvector of $G$ with eigenvalue $\alpha(G)$ and norm $\|\v\|=1$, so that
		\begin{equation}
			\label{lambda(v)=alpha}
			\lambda(\v)=\sum\limits_{(i,j)\in E}[\v(i)-\v(j)]^2=\alpha(G).
		\end{equation}
		Let $\w\in \Re^{n+1}$ be the vector with components $\w(i)=\v(i)-\frac{\v(n)}{n+1}$ for $i=1,\ldots,n$ and $\w(n+1)=\w(n)=\frac{n\v(n)}{n+1}$. Note that $\w$ is an embedding of $G^+$ since
		$$
		\sum_{i=1}^{n+1}\w(i)=\sum_{i=1}^n \left(\v(i)-\frac{\v(n)}{n+1}\right)+\w(n+1)=1-\frac{n\v(n)}{n+1}+\w(n)=1.
		$$
		Then, from \eqref{bound-lambda(v)},
		\begin{align*}
			\alpha(G^+) & \le \lambda(\w)=\frac{\sum\limits_{(i,j)\in E^+}[\w(i)-\w(j)]^2}{\sum\limits_{i\in V^+}\w^2(i)}
			=\frac{\sum\limits_{(i,j)\in E}[\v(i)-\v(j)]^2}{\sum\limits_{i\in V}[\v(i)-\frac{\v_n}{n+1}]^2+\w_{n+1}^2}\le \alpha(G)
		\end{align*}
		where the last inequality comes from \eqref{lambda(v)=alpha}  since, as $\v$ is an embedding of $G$,
		$$
		\sum\limits_{i\in V}\left[\v(i)-\frac{\v(n)}{n+1}\right]^2=\sum\limits_{i\in V}\left[\v(i)^2-2\v(i)\frac{\v(n)}{n+1}+\frac{\v(n)^2}{(n+1)^2}\right]=1+\frac{\v(n)^2}{(n+1)^2}\ge 1.
		$$
		Finally, the equality $\alpha(G^+)=\alpha(G)$ holds if and only if $\v(n)=0$.
	\end{proof}

	Let $G$ be a graph with $k$-token graph $F_k(G)$.
	For a vertex $a\in V(G)$, let $S_a:=\{A\in V(F_k(G)):a\in A\}$ and $S'_a:=\{B\in V( F_k(G)):a\not\in  B\}$.
	Let $H_a$ and $H'_a$ be the subgraphs of $F_k(G)$ induced by $S_a$ and $S'_a$, respectively.
	Note that $H_a\cong F_{k-1}(G\setminus \{a\})$ and $H'_a\cong F_k(G\setminus \{a\})$.
	\begin{lemma}
		\label{lem:eigenvectors}
		Given a vertex $a\in G$ and an eigenvector $\v$ of $F_k(G)$ such that $\B^{\top}\v=\vec0$, let
		\[
		\w_a:=\restr{\v}{S_a} \text{ and \quad } \w'_a:=\restr{\v}{S'_a}.
		\]
		Then, $\w_a$ and $\w'_a$ are embeddings of $H_a$ and $H'_a$, respectively.
	\end{lemma}
	
	\begin{proof}
		Assume that the matrix $\B^{\top}$ has the first row indexed by   $a\in V(G)$. Then, we have
		$$
		\vec0=
		\B^{\top}\v=
		\left(
		\begin{array}{cc}
			\1^{\top} & \vec0^{\top}\\
			\B_1 & \B_2
		\end{array}
		\right)
		\left(
		\begin{array}{c}
			\w_a \\
			\w'_a
		\end{array}
		\right)=
		\left(
		\begin{array}{c}
			\1^{\top}\w_a \\
			\B_1\w_a+\B_2\w'_a
		\end{array}
		\right),
		$$
		where $\1^{\top}$ is a row ${n-1\choose k-1}$-vector, $\vec0$ is a row  ${n-1\choose k}$-vector, $\B_1=\B(n-1,k-1)^{\top}$, and  $\B_2=\B(n-1,k)^{\top}$.
		Then, $\1^{\top}\w_a=0$, so that $\w_a$ is an embedding of $H_a$. Furthermore, since $\v$ is an embedding of $G$, we have $\1^{\top}\v=\1^{\top}\w_a+\1^{\top}\w'_a=0$ (with the appropriate dimensions of the all-1 vectors). Hence, it must be  $\1^{\top}\w'_a =0$, and $\w'_a$ is an embedding of $H'_a$.
	\end{proof}
	
	Now, we introduce some new results related to Conjecture \ref{conjecture}.

	\begin{theorem}
		\label{theo:alg-connec}
		For each  of the following classes of graphs, the algebraic connectivity of a token graph $F_k(G)$ satisfies the following.
		\begin{itemize}
			\item[$(i)$]
			Let $T_n$ be a tree on $n$ vertices. Then,
			$\alpha(F_k(T_n))=\alpha(T_n)$ for every $n$ and $k=1,\ldots,n-1$.
			\item[$(ii)$]
			Let $G$ be a graph such that $\alpha(F_k(G))=\alpha(G)$. Let $T_G$ be a graph where each vertex of $G$ is the root vertex of some (possibly empty) tree. Then
			$\alpha(F_k(T_G))=\alpha(T_G)$.
		\end{itemize}
	\end{theorem}
	
	\begin{proof}
		To prove $(i)$, let $V(T_n)=[n]$. From previous comments, we can assume that $T_n$ is connected. Then, the result is readily checked for $n\le 4$ and $k=1,2$. Now, we proceed by induction.
		Suppose $n>4$ and $k>1$. To our aim, by Lemma \ref{coro:LkL1}$(ii)$, it suffices to show that if $\v$ with a given norm, say $\v^{\top}\v=1$, is an eigenvector of $F_k:=F_k(T_n)$, with $\B^{\top}\v=\vec0$, then $\lambda(\v)\geq \alpha(T_n)$.
		Let $i\in [n]$. As defined before, let $S_i:=\{A\in V(F_k):i\in A \}$ and $S'_i:=\{B\in V(F_k):i\not\in B\}$.
		Let $H_i$ and $H'_i$ be the subgraphs of $F_k$ induced by $S_i$ and $S'_i$, respectively. We have $H_i\cong F_{k-1}(T_{n-1})$ and $H'_i\cong F_k(T_{n-1})$, where $T_{n-1}=T\setminus i$. Moreover, note that if vertex $i$ is of degree $1$ in $T_n$, then $T_{n-1}$ is also connected. Let $\w_i:=\restr{\v}{S_i}$ and $\w'_i:=\restr{\v}{S'_i}$, by Lemma~\ref{lem:eigenvectors},
		we know that $\w_i$ and $\w'_i$ are embeddings of $H_i$ and $H'_i$, respectively. By the induction hypothesis, we have
		\[
		\lambda(\w_i)=\frac{\sum\limits_{(A,B)\in E(H_i)} [\w_i(A)-\w_i(B)]^2}{\sum\limits_{A\in V(H_i)}\w_i(A)^2}\geq \alpha(T_{n-1}),
		\]
		and
		\[
		\lambda(\w'_i)=\frac{\sum\limits_{(A,B)\in E(H'_i)}[\w'_i(A)-\w'_i(B)]^2}{\sum\limits_{A\in V(H'_i)}\w'_i(A)^2}\geq \alpha(T_{n-1}).
		\]
		Since $V(H_i)\cup V(H'_i)=V(F_k)$ and $\v^{\top}\v=1$, we have
		\begin{align}
			\lambda(\v)&=\sum\limits_{(A,B)\in E(F_k)}[\v(A)-\v(B)]^2  \nonumber\\
			& \geq \sum\limits_{(A,B)\in E(H_i)}[\w_i(A)-\w_i(B)]^2 + \sum\limits_{(A,B)\in E(H'_i)}[\w'_i(A)-\w'_i(B)]^2 \nonumber\\
			& \geq \alpha(T_{n-1})\Big[\sum\limits_{A\in V(H_i)}\w_i(A)^2 + \sum\limits_{B\in V(H'_i)}\w'_i(B)^2\Big] \nonumber \\
			& = \alpha(T_{n-1})\Big[ \sum\limits_{A\in V(H_i)}\v(A)^2 + \sum\limits_{B\in V(H'_i)}\v(B)^2\Big] \nonumber \\
			& = \alpha(T_{n-1}) > \alpha(T_n), \label{eq:paths-1}
		\end{align}
		where (\ref{eq:paths-1}) follows  from Lemma \ref{lem:-vertex}. (Notice that, since $i$ has degree $1$, collapsing the edge of which $i$ is an end-vertex is equivalent to removing $i$, so obtaining $T_{n-1}$.)
		Furthermore, since $\lambda(\v)>\alpha(T_n)$, we get that $\alpha(T_n)$ is an eigenvalue of
		both $T_n$ and $F_k(T_n)$ with the same multiplicity.
		
		Regarding $(ii)$, it could be seen as a generalization of $(i)$. Thus, it is proved in the same way by induction on the number of vertices not in $G$ (that is, the non-root vertices of the trees), and starting from $G$.
		(The other way around, proved $(ii)$, the result in $(i)$ is a corollary when we start with $G=K_1$ or $G=K_2$.)
	\end{proof}
	
	The last step in (\ref{eq:paths-1}) also can be seen as a consequence of the following theorem
	by  Patra and Lal \cite[Th. 3.1]{pl08}.
	\begin{theorem}[\cite{pl08}]
		\label{th:pl08}
		Let $e=(u, v)$ be an edge of a tree $T$. Let $\widetilde{T}$ be the tree obtained from $T$ by `collapsing' the edge $e$ (that is, deleting $e$ and identifying $u$ and $v$). Then $\alpha(\widetilde{T})\ge \alpha(T)$.
	\end{theorem}
	
	Note 
	that the result of Theorem \ref{theo:alg-connec}$(i)$ implies the ones of Theorem \ref{theo:alg-connec-antic}$(ii)$ and $(iii)$.
	
	\begin{theorem}
		\label{theo:Delta}
		Let $G$ be a graph on $n$ vertices satisfying $\alpha(F_{k-1}(G))=\alpha(G)$ and maximum degree
		\begin{equation}
			\label{boundDelta}
			\Delta(G)\ge \phi(k)=\frac{k(n+k-3)}{2k-1}
		\end{equation}
		for some integer $k=1,\ldots, \lfloor n/2\rfloor$. Then, the algebraic connectivity of its $k$-token graph equals the algebraic connectivity of $G$,
		$$
		\alpha(F_k(G))=\alpha(G).
		$$
	\end{theorem}
	\begin{proof}
		When we `go' from the spectra of $\{F_{k-1}(G)$, $F_{k-1}(\overline{G})\}$ to the spectra of $\{F_{k}(G)$, $F_{k}(\overline{G})\}$,
		all the eigenvalues of $\Lambda_0,\ldots,\Lambda_{k-1}$ `reappear' (with eigenvectors $\v$ such that $\B^{\top}\v\neq \vec0$),
		together with `new' eigenvalues belonging to $\Lambda_k$ (with eigenvectors $\v$ such that $\B^{\top}\v= \vec0$). Moreover, the hypothesis $\alpha(F_{k-1}(G))=\alpha(G)$ implies that, in $\spec F_{k}(G)$, all
		the eigenvalues of $F_k(G)$ that are in $\Lambda_1,\ldots,\Lambda_{k-1}$ must be greater than or equal to $\alpha(G)$.
		Reasoning by contradiction, if $\alpha(F_k(G))<\alpha(G)$, then the eigenvalue $\alpha(F_k(G))$ must belong to 
		$\Lambda_k$. Then, the eigenvalue $\lambda_{F_k(G)}=\alpha(F_k(G))$ must be paired with one eigenvalue $\lambda_{F_k(\overline{G})}$ of 
		$F_k(\overline{G})$ belonging also to $\Lambda_k$ (both eigenvalues sharing the same eigenvector $\v$ with $J(n,k)$), so that 
		$$
		\alpha(G)+\lambda_{F_k(\overline{G})}>\alpha(F_k(G))+\lambda_{F_k(\overline{G})}=k(n-k+1)
		$$ 
		Thus,
		using that $\alpha(G)=n-\lambda_{\max}(\overline{G})$,
		$$
		\lambda_{\max}(F_k(\overline{G}))\ge \lambda_{F_k(\overline{G})}> k(n-k+1)-\alpha(G)=k(n-k+1)-n+\lambda_{\max}(\overline{G}).
		$$
		However, from the upper and lower bounds in \eqref{bound-sp-rad} for the spectral radius of a graph, together with \eqref{DeltaFk}, we get
		$$
		2k\Delta(\overline{G})\ge \lambda_{\max}(F_k(\overline{G}))>k(n-k+1)-n+\lambda_{\max}(\overline{G})\ge (k-1)(n-k)+\Delta(\overline{G})+1,
		$$
		or, in terms of $\Delta(G)$,
		$$
		n-1-\Delta(G)=\Delta(\overline{G})>\frac{(k-1)(n-k)+1}{2k-1}.
		$$
		Hence, $\Delta(G)<n-1-\frac{(k-1)(n-k)+1}{2k-1}=\frac{k(n+k-3)}{2k-1}$, contradicting the hypothesis.
	\end{proof}
	
	For the two extreme cases $k=2$ and $k=n/2$, we get the following consequences.
	
	\begin{corollary}
		Let $G$ be a graph on $n$ vertices and maximum degree $\Delta(G)$.
		\begin{itemize}
			\item[$(i)$]
			If $\Delta(G)\ge \frac{2}{3}(n-1)$, then $\alpha(F_2(G))=\alpha(G)$.
			\item[$(ii)$]
			If $\Delta(G)\ge \frac{3}{4}n$, then $G$ satisfies the Conjecture \ref{conjecture}. That is,
			$\alpha(F_k(G))=\alpha(G)$ for every $k=1,\ldots,n-1$.
		\end{itemize}
	\end{corollary}
	\begin{proof}
		$(i)$ With $k=2$, the condition \eqref{boundDelta} becomes $\Delta(G)\ge \frac{3}{2}(n-1)$. Then, since $\alpha(F_1(G))=\alpha(G)$, Theorem \ref{theo:Delta} gives the result.
		
		$(ii)$ Assuming that $n$ is even (the odd case is similar), it is enough to prove the result for $k=n/2$. In this case, the condition \eqref{boundDelta} becomes
		$\Delta(G)\ge \phi(n/2)=\frac{n(3n-6)}{4(n-1)}$. It is readily checked that $\frac{3}{4}n>\phi(n/2)>\phi(k)$ for every $k=2,\ldots,\frac{n}{2}-1$. So, we can use induction from the 
		case $(i)$ to prove the hypotheses in Theorem \ref{theo:Delta} hold for every $k$.
	\end{proof}
	
	Some examples of known graphs satisfying Conjecture \ref{conjecture} are:
	\begin{itemize}
		\item
		With maximum degree $n-1$, the complete and the star graphs (already mentioned), and the wheel graphs.
		\item
		With degree $n-2$, the cocktail party graph (obtained from the complete graph with even number of vertices minus a matching).
		\item
		With degree $n-3$, the complement $\overline{C_n}$ of the cycle with $n\ge 12$ vertices.
		\item
		The complete $r$-partite graph $G=K_{n_1,n_2,\ldots,n_r}\ne K_r$ for $r>2$, with number of vertices $n=n_1+n_2+\cdots +n_r$, for $n_1\le n_2\le \cdots \le n_r$ and $n\ge 4n_r$.
	\end{itemize}

\section*{Acknowledgments}
The authors are grateful to Clemens Huemer for his valuable comments.


\end{document}